\documentclass[11pt]{article}
\usepackage[utf8]{inputenc}
\usepackage[T1]{fontenc}
\usepackage{amsmath}
\usepackage{amsfonts}
\usepackage{amssymb}
\usepackage[version=4]{mhchem}
\usepackage{stmaryrd}
\usepackage{bbold}
\usepackage{hyperref}
\hypersetup{colorlinks=true, linkcolor=blue, filecolor=magenta, urlcolor=cyan,}
\numberwithin{equation}{section}
\newtheorem{Theorem}{Theorem}[section]

\newtheorem{Proposition}[Theorem]{Proposition}

\newtheorem{lemma}[Theorem]{Lemma}

\newenvironment{Proof}{\noindent{\bf Proof.}}{\hfill{$\blacksquare$}}

\title{Clark-Ocone formula for the maximum of processes with the stochastic intensity and its application}

\author{Mahdieh Tahmasebi}
\date{Department of Applied Mathematics, Tarbiat Modares University, P.O. Box 14115-175, Tehran, Iran.\\
\href{tahmasebi@modares.ac.ir}{tahmasebi@modares.ac.ir}
}

\begin{document}
\maketitle

\begin{abstract}
Pricing of the lookback options using the Clark-Ocone formula for the underlying assets driven by stochastic Lévy processes needs to compute the Malliavin derivatives of their maximum or minimum on the Wiener-Poisson space, as well as their distributions. In this work, we will 
find a generalization of explicit representation of the Clark-Ocone formula on the maximum of two types of Lévy processes with stochastic intensity; Cox processes with CIR-modeled intensities, and the Hawkes processes.   
\end{abstract}

\section*{1. Introduction and Preliminary Results}
The martingale representation theorem, or widely the Clark-Ocone formula in the context of Malliavin derivatives of Wiener-Poisson functionals, plays an important role in mathematical finance, such as in pricing of path-dependent options like lookback options, and in hedging a portfolio in risk management analysis, and so on.  For instance, in \cite{privault}, Privault and Teng established hedging strategies in bond markets, such as swaptions, using the Clark-Ocone formula with an suitable choice of numeraire. In 2017, Suzuki \cite{suzuki} obtained an explicit representation of a locally risk-minimizing hedging strategy in an incomplete financial market driven by multidimensional Lévy processes. Arai and Suzuki in 2021 \cite{arai} extended the Clark-Ocone theorem to random variables that are not Malliavin differentiable. They are given an explicit representation of the locally risk-minimizing strategy of digital options, particularly in the context of a 1-dimensional exponential Lévy process.  \\
Consider the completed probability space $(\Omega, \mathcal{F}, (\mathcal{F}_{t})_{t \in[0, T]}, \mathbb{P})$  satisfies $\mathcal{F}_{t-}=\mathcal{F}_{t}$, for every fixed time $t$, and let $\mathcal{P}$ be the predictable $\sigma$-field on $\mathbb{R} \times \Omega$ and $\mathcal{B}(\mathbb{R})$ the Borel $\sigma$-field on $\mathbb{R}$, for strictly positive real number $T$.
Denote $\mathcal{N}$ the space of simple and locally finite counting measures on $\mathbb{R} \times \Omega$, endowed with the vague topology, and $\nu_0$ Lebesgue measure on $[0, T]$. A process $\psi(t, z, \omega)$ is said to be Borel predictable if it is $( \mathcal{B}({[0, T]}) \times \mathcal{P} )$-measurable. Let $\bar{N}$ is a Poisson random measure with the Lévy finite maesure $\nu_1$ such that $\nu_1(\{0\})=0$, and 
the compensated random measure $\tilde{\bar{N}}$ is defined by

$$
\tilde{\bar{N}}([0, t] \times A)=\bar{N}([0, t] \times A)-t \nu_1(A).
$$
The following martingale representation theorem for Lévy processes demonstrated in \cite{tankov} Proposition 9.4.

\begin{Proposition} 
Let $F \in L^{2}(\Omega, \mathcal{F}, \mathbb{P})$. There exists a unique Borel predictable process $\psi \in L^{2}([0, T] \times \mathbb{R} \times \Omega)$ and a unique predictable process $\phi \in L^{2}([0, T] \times \Omega)$ such that
\begin{equation*}
F=\mathbb{E}[F]+\int_{0}^{T} \phi(t) d W_{t}+\int_{0}^{T} \int_{\mathbb{R}} \psi(t, z) \tilde{\bar{N}}(d t, d z). \tag{1.2}
\end{equation*}
\end{Proposition}
An extension of this representation can be expressed as the Clark-Ocone formula, based on Malliavin derivatives, Theorem 12.20 in \cite{nunno}. 
\begin{Proposition} 
For every random variale $F \in \mathbb{D}^{1,2}(\Omega, \mathcal{F}, \mathbb{P})$,
$$
F=\mathbb{E}[F]+\int_{0}^{T} \mathbb{E}\left[D_{t}^{(1)} F \mid \mathcal{F}_{t}\right] d W_{t}+\int_{0}^{T} \int_{\mathbb{R}} \mathbb{E}\left[D_{t, z}^{(2)} F \mid \mathcal{F}_{t}\right] \widetilde{\bar{N}}(d t, d z).
$$
\end{Proposition}
Here, two directional derivative operators $D^{(1)}: \mathbb{D}^{(1)} \rightarrow$ $L^{2}([0, T] \times \Omega)$ in the direction of the Brownian motion and $D^{(2)}: \mathbb{D}^{(2)} \rightarrow L^{2}([0, T] \times \mathbb{R} \times \Omega)$ in the direction of the Poisson random measure, where $\mathbb{D}^{(1)}$ and $\mathbb{D}^{(2)}$ stand for their respective domain,  are respectively given in the following next section. 
We say that $F$ is Malliavin differentiable if $F \in$ $\mathbb{D}^{1,2}:=\mathbb{D}^{(1)} \cap \mathbb{D}^{(2)}$. \\
For $0 \leq s<t \leq T$, define $M_{s, t}=\sup _{s \leq r \leq t} X_{r}$ and $M_{t}=M_{0, t}$. 
In \cite{13} have shown that if $X$ is a standard Brownian motion, then
$$
M_{T}=\sqrt{\frac{2 T}{\pi}}+\int_{0}^{T} 2\left[1-\Phi\left(\frac{M_{t}-W_{t}}{\sqrt{T-t}}\right)\right] d W_{t},
$$
where $\Phi(x)=\mathbb{P}\{N(0,1) \leq x\}$.    
The similar results hold for a Brownian motion with drift (see \cite{11}). For jump-diffusion Lévy processes, such as the Kou model \cite{5}, the inversion of a Laplace transform can obtain the distribution of $M_{T}$. Rémillard and  Renaud in 2011 \cite{remillard} computed an explicit representation for the maximum running of the assets driven by a square-integrable Lévy process in the form 
$$X_t= \mu t+\sigma W_t+ \int_0^t \int_{\mathbb{R}} z \tilde{\bar{N}}(ds,dz),$$
 where $\int_{\mathbb{R}} (1\wedge z^2) \nu(dz) < \infty$.\\
This setup is essentially used for pricing or hedging the lookback options in mathematical finance. 
Pricing and hedging a lookback option as a path-dependent option are considered by many authors. In Wiener space, we refer to \cite{bermin} and \cite{gobet}, in the fractional case we refer to \cite{Liu}.  
Also, Navarro in her thesis \cite{navarro}, using the Clark-Ocone formula, evaluated the Greeks of Exotic options with payoff of supremum or infimum type on the underlying assets driven by Lévy processes.     \\
Our main result is a generalization of the Clark-Ocone formula for point processes with stochastic intensity and provides an explicit expression for the Malliavin derivatives of lookback options, whee the underlying assets are driven by a Cox process or Hawkes process.\\ 
Cox processes and Hawkes processes are point processes with different behaviour of their stochastic intensity at jump times, describing greater flexibility in capturing the randomness of event occurrences. In a continuous Markov model, \cite{Damiano} derived an analytical formula for pricing credit derivatives under Cox-Ingersoll-Ross (CIR) stochastic intensity models. Subsequently, a smile-adjusted jump stochastic intensity to price credit default swaptions is extended by  \cite{brigo}. Additionally, non-Gaussian intensity models were developed in \cite{bianchi}. Moreover, the dynamics of risky asset prices have also  been modeled by jumps with self-exciting features in \cite{Filimonov et al. (2014)}. Such models produce large amount of price sudden movements triggered by previous sudden movements. Hawkes jump-diffusion models generate heavier-tailed distributions with higher peaks than Poisson jump-diffusion models, resulting in higher option prices for deep out-of-the-money options. 
Brignone and Sgarra \cite{brignone2020asian} have presented a method for pricing Asian options within risk asset models driven by Hawkes processes with exponential kernels.
In \cite{chen2024modeling}, the authors examine the clustering behavior of price jumps and their variances using high-frequency data modeled by a Hawkes process. We are interested in pricing lookback options in these models.\\
The rest of the paper is organized as follows. In Section 2, we recall some results from Malliavin calculus for Lévy processes.  Then, in Section 3, we state our main results on Cox processes and the distribution of the maximum running time of these processes by the inversion of a Laplace transform. In Section 4, we present an explicit formula for the maximum running of the Hawkes processes and present their distribution by the inversion of a Laplace transform.


\section{Malliavin calculus on Wiener-Poisson Space}
In this section, we recall the concept of Malliavin calculus from \cite{nualart}. Consider the Wiener-Poisson space $(\Omega, \mathcal{F}, \mathcal{F}_t, P)$. 
Given $\gamma(t)=\int_{0}^{t}g(s)ds$ for some $g\in{L^{2}}(\left[0, T\right])$, and a random variable $F:\Omega\to{\mathbb{R}}$, the directional derivative of $F$ in the $\gamma$ direction in Wiener space, have defined as the following form, if the limit exists,  e.g.,
	\begin{equation*}
	D^{(1)}_\gamma{F}(w)=\frac{d}{d\epsilon}[F(w+\epsilon\gamma)]_{\epsilon=0}.
	\end{equation*}
The variable $F$ is Malliavin differentiable in Wiener space, if there exists some $\psi\in{L^{2}(\left[0,T\right]\times\Omega)}$ satisfying the following equation
	\begin{equation}
	\nonumber D^{(1)}_\gamma{F}(w)=\int_{0}^{T}\psi(t,w).g(t)dt.
	\end{equation}
Denote $D^{(1)} F:=(\psi(t,w))_{0\le{t}\le{T}}$. We define the set of all $F:\Omega\to\mathbb{R}$ such that $F$ is differentiable in Wiener space by $\mathbb{D}^{(1)}$ with the norm defined by 
\begin{equation*}\label{equ20}
\lVert{F}\rVert^2=\lVert{F}\rVert^2_{L^{2}(\Omega)}+\lVert{D^{(1)}{F}}\rVert^2_{L^{2}(\left[0, T\right]\times\Omega)}.
\end{equation*}
 Given $h\in L^{2}([0,T] \times \mathbb{R}_0^{n})$ and fixed $z\in \mathbb{R}_0$, we write $h(t,.,z)$ to indicate the function on $\mathbb{R}_0^{n-1}$ given by $(z_1,...,z_{n-1})\to h(t, z_1,...,z_{n-1}, z)$.
Denote by $\mathbb{D}^{(2)}$ the set of random variables $F$ in $L^{2}(\Omega)$ with a chaotic decomposition $F=\sum_{n=0}^{\infty}I_n(h_n)$, that $h_n\in L_s^{2}([0,T] \times \mathbb{R}_0^{n})$ and $I_n$ is a $n$-dimensional Poisson integral, satisfying
\begin{equation*}
	\nonumber \sum_{n\ge1}nn!\lVert{h_n}\rVert_{L^{2}([0,T] \times \mathbb{R}_0^{n})}^{2}<\infty.
\end{equation*}
For every $F\in \mathbb{D}^{(2)}$, the Malliavin derivative of $F$ in a Poisson space is defined as the $L^{2}([0,T] \times \mathbb{R}_0)$-valued random variable given by
\begin{equation*}
	\nonumber D_{t,z}^{(2)}F=\sum_{n\ge1}nI_{n-1}(h_n(t, .,z)),\;\;z\in \mathbb{R}_0.
\end{equation*}
Let $\mathbb{D}^{1,2}= \mathbb{D}^{(1)} \cap \mathbb{D}^{(2)}$, and for the Malliavin derivative $DF:=\left(D^{(1)} F, D^{(2)} F\right)$, put the norm 
$$
\|D F\|^{2}=\left\|D^{(1)} F\right\|_{L^{2}([0, T] \times \Omega)}^{2}+\left\|D^{(2)} F\right\|_{L^{2}([0, T] \times \mathbb{R}_0 \times \Omega)}^{2}.
$$

The Malliavin derivative $D$ is a closed operator.
\begin{Proposition}\label{closed}\cite{nualart}
 If $F$ belongs to $L^{2}(\Omega)$, if $\left(F_{k}\right)_{k \geq 1}$ is a sequence of elements in $\mathbb{D}^{1,2}$ converging to $F$ in the $L^{2}(\Omega)$-norm and if $\sup _{k \geq 1}\left\|D F_{k}\right\|<\infty$, then $F$ belongs to $\mathbb{D}^{1,2}$ and $\left(D F_{k}\right)_{k \geq 1}$ converges weakly to $D F$ in $L^{2}([0,T] \times \Omega) \times L^{2}([0,T] \times \mathbb{R}_0 \times \Omega)$.
\end{Proposition}

\begin{Proposition}\label{chainf}\cite{nualart}
	Let $F$ be a random variable in $\mathbb{D}^{1,2}$ and let $\varphi$ be a real continuous function such that $\varphi(F)$ belongs to $L^{2}(\Omega)$ and $\varphi(F + D^{N}F)$ belongs to $L^{2}([0,T] \times \mathbb{R}_0 \times \Omega)$. Then $\varphi(F)$ belongs to $\mathbb{D}^{1,2}$ and
	\begin{equation}
		\nonumber D_t^{(1)}\varphi(F)=\varphi'(F) D_{t}^{(1)}F, \qquad \nonumber D_{t,z}^{(2)}\varphi(F)=\varphi(F+D_{t,z}^{(2)}F)-\varphi(F).
	\end{equation}
\end{Proposition}
As a consequence, if $G=g\left(X_{t_{1}}, \ldots, X_{t_{n}}\right) \in \mathbb{D}^{(2)}$ for some Lipschitz function $g$  on $\mathbb{R}^n$, and 
$$
(t, z) \mapsto g\left(X_{t_{1}}+D^{(2)}_{t,z}X_{t_1}, \ldots, X_{t_{n}}+D^{(2)}_{t,z}X_{t_n}\right)-g\left(X_{t_{1}}, \ldots, X_{t_{n}}\right),
$$
belongs to $L^{2}([0,T] \times \mathbb{R}_0 \times \Omega)$, then
$$
D_{t, z}^{(2)} G=g\left(X_{t_{1}}+D^{(2)}_{t,z}X_{t_1}, \ldots, X_{t_{n}}+D^{(2)}_{t,z}X_{t_n}\right)-g\left(X_{t_{1}}, \ldots, X_{t_{n}}\right).
$$
In addition, let $\left(t_{k}\right)_{k \geq 1}$ be a dense subset of $[0, T]$, $F=M_{T}=\max_{0 \le t \leq T}X_t$, and for each $n \geq$ 1, define $F_{n}=\max \left\{X_{t_{1}}, \ldots, X_{t_{n}}\right\}$. Clearly, $\left(F_{n}\right)_{n>1}$ is an increasing sequence bounded by $F$. Hence $F_{n}$ converges to $F$ in the $L^{2}(\Omega)$-norm when $n$ goes to infinity. If the process $X_t$ belongs to $\mathbb{D}^{1,2}=\mathbb{D}^{(1)} \cap \mathbb{D}^{(2)}$, since
$$
\left(x_{1}, \ldots, x_{n}\right) \mapsto \max \left\{x_{1}, \ldots, x_{n}\right\}
$$
is a Lipschitz function on $\mathbb{R}^{n}$, from Proposition \ref{closed}, we know that $DF_n$  belongs to $\mathbb{D}^{1,2}$ and finally if the uniformly boundedness of $DF_n$ holds, due to Proposition \ref{chainf}, the Malliavin differentiability of $F$ will be concluded. 
\\
%
Now, we are ready to show an extension of the martingale representation for the supremum of two types of processes: Cox processes and Hawkes processes, in the next sections. 
\section{Cox processes}
Recall the concept of stochastic intensity desired by B{\'e}rmaud in Chapter 5 of \cite{bremaud}. 
For given  $\sigma$-field $\mathcal{G}_t$, the process $\lambda_t$ is an $\mathcal{G}_t$-intensity of a Poisson process $N^c_t$ if for every $s,t \in [0,T]$
\begin{equation*} 
	\mathbb{E}\Big(\int_{\mathbb{R}_0}\int_t^s  N^c(du,dz)\vert \mathcal{G}_t\Big)=\mathbb{E}\Big(\int_{\mathbb{R}_0}\int_t^s \lambda_u du\nu(dz)\vert \mathcal{G}_t\Big),
\end{equation*}
and so that  $\tilde{N}^c(t,z)=N^c(t.z)-\int_{\mathbb{R}_0}\int_0^t\lambda_sds\nu(dz)$  is an $\mathcal{G}_t$-martingale. Also, for every $0 \leq t,s \leq T$ and for every $\mathcal{G}_t$-predictable function $k$  
\begin{equation*} 
	\mathbb{E}\Big(\int_{\mathbb{R}_0}\int_t^s k(u,z) N(du,dz)\vert  \mathcal{G}_t \Big)=\mathbb{E}\Big(\int_{\mathbb{R}_0}\int_t^s k(u,z)\lambda_u du\nu(dz)\vert  \mathcal{G}_t \Big). 
\end{equation*}
In this section, we suppose that 

\begin{equation}\label{Slambda1}
\begin{cases}
dS_{t}&=\mu  S_{t}dt+\sigma_{1}S_{t}dW_{t}^S+\int_{\mathbb{R}_{0}}(e^{J_{t,z}}-1)S_{t}\tilde{N}^c(dt,dz), \\
d\lambda_{t}&=\kappa(\Theta-\lambda_{t})dt+\sigma_{2}\sqrt{\lambda_{t}}dW_{t},
\end{cases}
\end{equation}
where $\mu$ is a real number, $\sigma_1, \sigma_2, \kappa$ and $\Theta$ are strictly positive real numbers, $W$ is a standard Brownian motion and is independent of the Brownian motion $W^S$ and $N^c$. \\
We assume that $$\bar{\mu}_t := \int_{\mathbb{R}_0} (e^{J_{t,z}}-1)\nu(dz) < \infty, \quad \int_0^T \int_{\mathbb{R}_{0}} \vert J_{t,z} \vert^2 \nu(d z)dt < \infty. $$
Denote  $X_{t}= ln S_t$ and rewrite the above equation of $S_t$ as the form 
\begin{equation*}
X_{t}=(\mu -\frac{\sigma^2}{2})t-\int_0^t \bar{\mu}_s\lambda_s ds +\sigma_1 W^S_{t}+\int_{0}^{t} \int_{\mathbb{R}} J_{s,z} {N}^c(d s, d z).
\end{equation*}
As mentioned in the previous section, for a dense subset $\left(t_{k}\right)_{k \geq 1}$  of $[0, T]$, and $F=M_{T}$, the processes $F_{n}=\max \left\{X_{t_{1}}, \ldots, X_{t_{n}}\right\}$ converge to $F$ in the $L^{2}(\Omega)$-norm when $n$ goes to infinity.
We want to prove that each $F_{n}$ is Malliavin differentiable and to find an expression of the supremum of the Cox process $S_t$ in \eqref{Slambda1} by using the Clark-Ocone formula. \\
Here, we mention that in \cite{altmayer}, the authors have shown that using the It\^o formula and taking the Malliavin derivative with respect to the Brownian motion $W$,  for every $s \leq t$
\begin{equation}\label{derlamb}
D_s^{(1)}\lambda_t=\sigma_{2}\sqrt{\lambda_t}1_{0\leq s\leq t}\exp\Big\{-{\int_{s}^{t}(\frac{\kappa}{2}+\frac{C_{\sigma}}{\lambda_r})dr}\Big\},
\end{equation}
where $C_{\sigma}=\frac{\kappa\Theta}{2}-\frac{\sigma_{2}^{2}}{8}$ is a positive number. 
Also, for all $p \geq 1$,
\begin{equation}\label{suplambda}
	\mathbb{E}(\sup_{0 \leq t \leq T} \lambda_t ^p) <\infty,  \quad    \quad  and   \quad  \sup_{0 \leq t \leq T}\mathbb{E}(\lambda_t ^{-p}) <\infty, \quad~if ~  2\kappa\Theta > p\sigma_2^2.
\end{equation}
From the definition of Malliavin derivative, we know that 
$$
0 \leq D_{t}^{(1)} F_{n}=\sum_{k=1}^{n} \int_t^{t_k} \bar{\mu}_s D_{t}^{(1)} \lambda_s ds \ \mathbb{I}_{\left\{t \leq t_{k}\right\}} \mathbb{I}_{A_{k}} \leq \sum_{k=1}^{n} \sigma_2 \int_t^{t_k} \bar{\mu}_s \sqrt{\lambda_s} ds \mathbb{I}_{A_{k}},
$$
where $A_{1}=\left\{F_{n}=X_{t_{1}}\right\}$, $A_{k}=\left\{F_{n} \neq X_{t_{1}}, \ldots, F_{n} \neq X_{t_{k-1}}, F_{n}=X_{t_{k}}\right\}$, for $2 \leq k \leq n$, and we used from (\ref{derlamb}) in the last inequality. This implies that 
\begin{align*}
\sup _{n \geq 1}\left\|D^{(1)} F_{n}\right\|_{L^{2}([0, T] \times \Omega)}^{2} & \leq \sigma^{2} \left\|\sup _{T \geq s \geq 0}  \sqrt{\lambda_s}\right\|_{L^{2}([0, T] \times \Omega)}^{2}  \sum_{k=1}^{n}\int_t^{t_k} \bar{\mu}_s  ds \mathbb{I}_{A_{k}} \\
&\leq \sigma_2^{2} \left\|\sup _{T \geq s \geq 0}  \sqrt{\lambda_s}\right\|_{L^{2}([0, T] \times \Omega)}^{2}\int_t^{T} \bar{\mu}_s  ds  .
\end{align*} 
Secondly, since $D^{(2)}$ operates as the Malliavin derivative on the Poisson part of $F_{n}$, we have that
\begin{align*}
D_{t, z}^{(2)} F_{n}&=\max \left\{X_{t_{i}}+\big(J_{t,z}-\int_t^{t_i} \bar{\mu}_s D^{(2)}\lambda_s ds\big)\mathbb{I}_{\left\{t<t_{i}\right\}}, i=1,\cdots k\right\}-F_{n}\\
& =\max \left\{X_{t_{i}}+J_{t,z}\mathbb{I}_{\left\{t<t_{i}\right\}}, i=1,\cdots k\right\}-F_{n},
\end{align*}
where the equality is justified by the following inequality:
\begin{align*}
&\left\|\max \left\{X_{t_{i}} +J_{t,z}\mathbb{I}_{\left\{t<t_{i}\right\}}, i=1,\cdots k\right\}-F_{n}\right\|_{L^{2}([0, T] \times \mathbb{R} \times \Omega)}^{2}\\ 
&\leq \int_0^T \int_{\mathbb{R}_{0}}\vert  J_{t,z}\vert^2 \nu(d z)dt < \infty.
\end{align*}
Consequently, 
\begin{align*}
\sup _{n \geq 1}\left\|D F_{n}\right\|^{2} &\leq  \sigma_2^{2} \left\|\sup _{T \geq s \geq 0}  \sqrt{\lambda_s}\right\|_{L^{2}([0, T] \times \Omega)}^{2}\int_t^{T} \bar{\mu}_s  ds \\
& +  \int_0^T \int_{\mathbb{R}_{0}} \vert J_{t,z} \vert^2 \nu(d z)dt < \infty,
\end{align*}
 and we have that $F$ is Malliavin differentiable. By the uniqueness of the limit, this means that taking the limit of $D_{t}^{(1)} F_{n}$ when $n$ goes to infinity yields $$D_{t}^{(1)} F= \mathbb{I}_{[0, \tau]}(t)\int_t^\tau \bar{\mu}_s D_t^{(1)}\lambda_s ds,$$
where $\tau=$ $\inf \left\{t \in[0, T]: X_{t}=M_{T}\right\}$, with the convention inf $\emptyset=T$, i.e. $\tau$ is the first time when the Lévy process $X$ reaches its supremum on $[0, T]$, and
\begin{equation}\label{d2cox}
D_{t, z}^{(2)} F=\sup _{0 \leq s \leq T}\left(X_{s}+J_{t,z} \mathbb{I}_{\{t<s\}}\right)-M_{T}.
\end{equation}
Now, we are ready to present the Clark-Ocone formula in the following theorem.  Set $\bar{F}_{t, x, y}=\mathbb{P}\left\{M_{t,T}>x, \sup_{t \leq s \leq T} \lambda_s > y\right\}$. 
\begin{Theorem}
For the square-integrable Cox process $S_t$, the running maximum of $log S_t$ can be written as follows:
\begin{align*}
M_{T}=\mathbb{E}\left[M_{T}\right] &-\int_{0}^{T} \mathbb{E}\left[ \mathbb{I}_{[0, \tau]}(t)\int_t^\tau \bar{\mu}_s D_t^{(1)}\lambda_s ds\mid \mathcal{F}_{t}\right] dW_t \\
& +\int_{0}^{T} \int_{\mathbb{R}} \psi\left(t, z, M_{t}-J_{t, z}\right) \tilde{N}^c(d t, d z),
\end{align*}
where $\psi(t, z, y)=\mathbb{I}_{\{J_{t,z} \geq 0\}}\int_{y}^{M_t} \int_{0}^{\infty}\bar{F}_{t, x, r}dr d x -  \mathbb{I}_{\{J_{t,z} < 0\}}\int_{M_t}^{y} \int_{0}^{\infty}\bar{F}_{t, x, r}dr d x$.
\end{Theorem}
\begin{Proof}
From the above discussion, we know 
$$
\mathbb{E}\left[D_{t}^{(1)} F \mid \mathcal{F}_{t}\right]= \mathbb{E}\left[\mathbb{I}_{[0, \tau]}(t)\int_t^\tau \bar{\mu}_s D_t^{(1)}\lambda_s ds\mid \mathcal{F}_{t}\right].
$$
Also, using Equation \eqref{d2cox} we get that
\begin{align*}
\mathbb{E}\left[D_{t, z}^{(2)} F \mid \mathcal{F}_{t}\right] & =\mathbb{E}\left[\sup _{0 \leq s \leq T}\left(X_{s}+J_{t,z} \mathbb{I}_{\{t<s\}}\right)-M_{T} \mid \mathcal{F}_{t}\right] \\
& =\mathbb{E}\left[\max \left\{M_{t}, M_{t, T}+J_{t,z}\right\} \mid \mathcal{F}_{t}\right]-\mathbb{E}\left[M_{T} \mid \mathcal{F}_{t}\right] \\
& =M_{t}+\mathbb{E}\left[\left(M_{t, T}+J_{t,z}-M_{t}\right)^{+} \mid \mathcal{F}_{t}\right]-\mathbb{E}\left[M_{T} \mid \mathcal{F}_{t}\right] \\
& =\mathbb{E}\left[\left(M_{t, T}-a\right)^{+} \mid \mathcal{F}_{t}\right]-\mathbb{E}\left[(M_{t, T} - M_t )^{+} \mid \mathcal{F}_{t}\right],
\end{align*}
where $a=M_{t}-J_{t,z}$. Since 
\begin{align*}
\mathbb{E}\left[\left(M_{t, T}-a\right)^{+} \mid \mathcal{F}_{t}\right]  &=\int_{a}^{\infty}  P(\sup_{t \leq s \leq T}X_s > x \mid  \mathcal{F}_{t})d x\\
& =\int_{a}^{\infty} \int_{0}^{\infty}  P(\sup_{t \leq s \leq T}X_s > x , \sup_{t \le s \le T}\lambda_s > y \mid  \mathcal{F}_{t})dy d x\\
&=: \int_{a}^{\infty}\int_{0}^{\infty} \bar{F}_{t, x, y}dy d x,
\end{align*}
we have that
\begin{align*}
\mathbb{E}\left[D_{t, z}^{(2)} F \mid \mathcal{F}_{t}\right] = & \int_{a}^{\infty}\int_{0}^{\infty} \bar{F}_{t, x, y}dy d x- \int_{M_t}^{\infty}\int_{0}^{\infty} \bar{F}_{t, x, y}dy d x\\
= & \mathbb{I}_{\{J_{t,z} \geq 0\}}\int_{a}^{M_t} \int_{0}^{\infty}\bar{F}_{t, x, y}dy d x -  \mathbb{I}_{\{J_{t,z} < 0\}}\int_{M_t}^{a} \int_{0}^{\infty}\bar{F}_{t, x, y}dy d x\\
=& \psi(t, z, a).
\end{align*}
\end{Proof}

\noindent Due to the similarity of the proof to obtain the distribution function $\bar{F}_{., ., .}$ in the Cox process with that in the Hawkes process, we will consider this part in Appendix.
\begin{lemma}\label{distcox1}
There exist some positive constants $\alpha_1$ and $\alpha_2$ such that 
$$\bar{F}_{t,x,y} =   P({\tau}_{t}^0(x,y) \leq T \vert \mathcal{F}_t ) =  \frac{1}{\alpha_2} L^{-1}(\frac{1}{u})(T-t)  e^{-\alpha_1 (x-X_t)  -\alpha_2(y-\lambda_t)},$$
where $L^{-1}(f)(x)$ is the inverse Laplace operator of the function $f$ in the point $x$.
\end{lemma}
\begin{Proof}
We will prove in Appendix.
\end{Proof}


\section{Hawkes process}
 Consider point processes $\mathit{N}$ on $\mathbb{R}$ defined by
\begin{equation*}
N(dt):=\overline{N}(dt\times\left(0,\lambda_t\right]),
\end{equation*}
where $\{\lambda_t\}_{t\in\mathbb{R}}$ is a nonnegative process of the form
\begin{equation*}
\lambda_t:=\varphi(t,\overline{N}\lvert_{\left(-\infty,t\right)})=\varphi_t(\overline{N}\lvert_{\left(-\infty,t\right)}),
\end{equation*}
such that for all $a,b\in\mathbb{R}$ we have $\int_{a}^{b}\lambda_t dt<\infty$ almost surely.
Here, the map $\varphi:\mathbb{R}\times\mathcal{N}\to\mathbb{R_{+}}$ is a measurable functional, $\mathcal{N}$ denotes the space of simple and locally finite counting measures on $\mathbb{R}\times\mathbb{R_{+}}$ endowed with the vague topology. For simplicity, with a little abuse of notation, we denote by $\overline{N}\lvert_{\left(-\infty,t\right)}$ the restriction of $\overline{\mathit{N}}$ to ${\left(-\infty,t\right)}\times\mathbb{R_{+}}$. 
Since the process $\{\overline{N}\lvert_{\left(-\infty,t\right)}\}_{t\in\mathbb{R}}$ is $\mathcal{F^{\overline{N}}}$-adapted and left-continuous, therefore
$\{\lambda_t\}_{t\in\mathbb{R}}$ is $\mathcal{F^{\overline{N}}}$-predictable. Consequently, Lemma 2.1 in \cite{torrisi} deduces that $N$ 
has $\mathcal{F^{\overline{N}}}$-stochastic intensity $\{\lambda_t\}_{t\in\mathbb{R}}$. Additionally, for every measurable function $u$, 
$$\int_{[0, T]} u(t) (N(dt) -\lambda_t dt) =\int_{[0, T]} \int_{\mathbb{R}_0} u(t)\mathbb{1}_{(0, \lambda_t]}(\bar{N}(dt,dz) - dtdz).$$ 
Consider the following stochastic differential equation  
 \begin{equation*}
\begin{cases} 
dS_{t}&=(\mu+(1-e^{J_t})\lambda_t)S_{t}dt+\sigma_{1}S_{t}dW_{t}+(e^{J_t}-1)S_{t}{N}(dt), \\
d\lambda_{t}&=\kappa(\Theta-\lambda_{t})dt+\eta N(dt).
\end{cases}
\end{equation*}
We assume that 
$$\mathbb{E}\Big(\int_0^T (e^{J_t}-1) \lambda_t dt\Big) < \infty, \quad \mathbb{E}\Big(\int_0^T (e^{J_t}-1)^2 \lambda_t dt\Big) < \infty.$$
It\^o formula drives that
\begin{align*}
\nonumber dY_t:= d (ln S_{t})=(\mu-\frac{1}{2}\sigma_1^{2})dt-(e^{J_t}-1)\lambda_t dt+\sigma_1 dW_{t}+J_ t N(dt).\\
\end{align*}
In \cite{torrisi}, the authors with using the It\^o formula and taking the Malliavin derivative on Wiener-Poisson space  have driven that for every $s \leq t$
\begin{equation}\label{14}
D_s^{(1)}\lambda_t=0, \qquad D_{s,z}^{(2)}\lambda_t =  \eta e^{\kappa(s-t)}\mathbb{I}_{(0, \lambda_s)}(z) +\int_s^t \eta e^{\kappa(u-t)} sgn(D_{s,z}^{(2)} \lambda_u) N_{s,z}(du),
\end{equation}
\begin{align}
\nonumber N_{(u,z)}(ds)=&\overline{N}(ds\times(\varphi_s(\overline{N}\lvert_{\left(-\infty,s\right)}+\epsilon_{(u,z)})\land\varphi_s(\overline{N}\lvert_{\left(-\infty,s\right)}),\\ \nonumber &\varphi_s(\overline{N}\lvert_{\left(-\infty,s\right)}+\epsilon_{(u,z)})\lor\varphi_s(\overline{N}\lvert_{\left(-\infty,s\right)})]),
\end{align}
where $\epsilon_{(u,z)}$ denotes the Dirac measure at $(u,z)\in\mathbb{R}\times\mathbb{R_+}$ and $a\land b$ and $a\lor b$ are the minimum and the maximum between $a,b\in\mathbb{R}$, respectively.
They also have shown that the $D\lambda$ is positive when the function $\varphi$ is an increasing function. 
 
Let $X_t=Y_t+\int_0^t (e^{J_u}-1) \lambda_u du$ and compute the Malliavin derivatives of $X_t$ in the Clark-Ocone formula for $F= M_T= max_{0 \leq s \leq T}X_s$. For every $t\leq T$
$$
0 \leq D_{t}^{(1)} F_{n}=\sum_{k=1}^{n} \sigma_1 \mathbb{I}_{\left\{t \leq t_{k}\right\}} \mathbb{I}_{A_{k}} \leq \sum_{k=1}^{n} \sigma_1 \mathbb{I}_{A_{k}}=\sigma_1,
$$
where $A_{1}=\left\{F_{n}=X_{t_{1}}\right\}$ and $A_{k}=\left\{F_{n} \neq X_{t_{1}}, \ldots, F_{n} \neq X_{t_{k-1}}, F_{n}=X_{t_{k}}\right\}$ for $2 \leq k \leq n$. This implies that 
\begin{align*}
\sup _{n \geq 1}\left\|D^{(1)} F_{n}\right\|_{L^{2}([0, T] \times \Omega)}^{2} & \leq \sigma_1^{2}.
\end{align*} 
Secondly, since $D^{(2)}$ operates like the Poisson random measure Malliavin derivative on the Poisson part of $F_{n}$, we have that
$$D_{t, z}^{(2)} X_{t_i}=J_{t}\mathbb{I}_{(0,\lambda_t)}(z)+\int_t^{t_i} J_s sign{(D_{t,z}^{(2)}\lambda_s)}N_{t,z}(ds), $$
$$
D_{t, z}^{(2)} F_{n}=\max \left\{X_{t_{i}}+D_{t,z}^{(2)}X_{t_i} \mathbb{I}_{\left\{t<t_{i}\right\}}, i=1,\cdots k\right\}-F_{n},
$$
where the equality is justified by the following inequality:
\begin{align*}
&\left\|\max \left\{X_{t_{i}} +D_{t,z}^{(2)}X_{t_i} \mathbb{I}_{\left\{t<t_{i}\right\}}, i=1,\cdots k\right\}-F_{n}\right\|_{L^{2}([0, T] \times \mathbb{R} \times \Omega)}^{2}\\ 
&\leq C_1T \sup_{t \leq s \leq T} \| D_{t,z}^{(2)}X_{s} \|_{L^{2}([0, T] \times \mathbb{R} \times \Omega)}^{2} < \infty.
\end{align*}
Consequently, $\sup _{n \geq 1}\left\|D F_{n}\right\|^{2} \leq T\left(\sigma_1^{2}+C_1\sup_{t \leq s \leq T} \| D_{t,z}^{(2)}X_s \|_{L^{2}([0, T] \times \mathbb{R} \times \Omega)}^{2}\right)$ and we have that $F$ is Malliavin differentiable. By the uniqueness of the limit, this means that taking the limit of $D_{t}^{(1)} F_{n}$ when $n$ goes to infinity yields $$D_{t}^{(1)} F=\sigma_1 \mathbb{I}_{[0, \tau]}(t),$$
where $\tau=$ $\inf \left\{t \in[0, T]: X_{t}=M_{T}\right\}$, with the convention inf $\emptyset=T$, i.e. $\tau$ is the first time when the Lévy process $X$ (not the Brownian motion $W$ ) reaches its supremum on $[0, T]$, and

\begin{align}
D_{t, z}^{(2)} F& =\sup _{0 \leq s \leq T}\left(X_{s}+J_{t}\mathbb{I}_{(0,\lambda_t)}(z)+\int_t^{s}J_{u}N_{t,z}(du)\right)-M_{T}\nonumber\\
&= \max \big\{M_t , \sup_{t \leq s \leq T}Z_s +J_t \mathbb{I}_{(0,\lambda_t)}(z)\big\}-M_{T}\nonumber\\
&=: \max \big\{M_t , \sup_{t \leq s \leq T}Z_s+ K_{t,z}\big\}-M_{T}, \label{d2hawkes}
\end{align}
where for every $t \leq s \leq T$
\begin{align*}
Z_s= X_t &+  (\mu - \frac12 \sigma^2) (s-t) +\sigma_1 (W_{s}-W_t) \\
&+\int_t^s \int_{\mathbb{R}-0} J_u \mathbb{I}_{(0, \varphi_u (\bar{N}|_{(-\infty, u)}+\epsilon_{t,z}))}(z_0)\bar{N}(du, dz_0).\\
\end{align*} 
Now, we are ready to present the Clark-Ocone formula in the following theorem.
\begin{Theorem}
For the square-integrable Hawkes process $S_t$, the running maximum of $log S_t$ can be written as follows:
\begin{align*}
M_{T}=\mathbb{E}\left[M_{T}\right] &-\int_{0}^{T} \sigma_1  P(M_t \leq M{_{t,T}}) dW_t \\
& +\int_{0}^{T} \int_{\mathbb{R}} \psi\left(t, z, M_{t}-K_{t, z}\right) \tilde{\bar{N}}(d t, d z),
\end{align*}
where $\psi(t, z, y)=\int_y^\infty P(\sup_{t \leq s \leq T}Z_s \geq y \vert \mathcal{F}_t) dy-\int_{M_t}^{\infty} P(M_{t,T} \geq x) d x$.
\end{Theorem}
\begin{Proof}
From above discussion, we know 
$$
\mathbb{E}\left[D_{t}^{(1)} F \mid \mathcal{F}_{t}\right]=\sigma_1 \mathbb{E}\left[\mathbb{I}_{[0, \tau]}(t)\mid \mathcal{F}_{t}\right]=\sigma_1  P(M_t \leq M{_{t,T}}).
$$
Also, using Equation \eqref{d2hawkes}, we get that
\begin{align}
\mathbb{E}\left[D_{t, z}^{(2)} F \mid \mathcal{F}_{t}\right] 
& =M_{t}+\mathbb{E}\left[\left(\sup_{t \leq s \leq T}Z_s+K_{t,z}-M_{t}\right)^{+} \mid \mathcal{F}_{t}\right]-\mathbb{E}\left[M_{T} \mid \mathcal{F}_{t}\right] \nonumber\\
& =\mathbb{E}\left[\left(\sup_{t \leq s \leq T}Z_s-a \right)^{+} \mid \mathcal{F}_{t}\right]-\int_{M_t}^{\infty} P(M_{t,T} \geq x) d x \nonumber\\
& =\int_a^\infty P(\sup_{t \leq s \leq T}Z_s \geq y \vert \mathcal{F}_t) dy-\int_{M_t}^{\infty} P(M_{t,T} \geq x) d x,  \label{distribute}
\end{align}
where $a=M_t - K_{t,z}$. 
\end{Proof}\\
\subsection{Distribution of the supremum running of Hawkes processes}
In this subsection, we find the distributions of $M_{t,T}$ and $\sup_{t \leq s \leq T}Z_s$ when the funcions $J$ and $\eta$ do not depend on the time, and $\eta < (\ln2)\kappa$. To this end, we first define the following two exit times.
Given constants $b$ and $e$, the random time $\hat{\tau}_{t,b,e}$ is the first time when the process $Z$ is more than $b$, the process $\lambda$ is more than $e$ and the process $D_{t,z}\lambda$ is more than $r$  on $[t,\infty]$.\\
Also, the random time ${\tau}_{t}(b,e)$ is the first time when the process $X$ is more than $b$ and the process $\lambda$ is more than $e$.\\
Let 
\begin{align*}
&\mathcal{L}_t^1 u(Z_s, \lambda_s, D_{t,z}^{(2)}\lambda_s)=-\kappa \lambda_s \partial_2 u(Z_s,\lambda_s, D_{t,z}^{(2)}\lambda_s) - \kappa \lambda_s \partial_3 u(Z_s,\lambda_s, D_{t,z}^{(2)}\lambda_s) \\
&+ \int_{\mathbb{R}_0}\Big(u(Z_s+\triangle Z_s, \lambda_s+\triangle \lambda_s, D_{t,z}^{(2)}\lambda_s+\triangle D_{t,z}^{(2)}\lambda_s)-u(Z_s,\lambda_s, D_{t,z}^{(2)}\lambda_s)\Big) dz_0,
\end{align*}
where $\partial_i u(x,y,w)$ is the partial deriative of $u$ with respect to $i$-th component, for $i=1,2,3$.
Put the positive constants $\alpha_1$ and $\alpha_2$ such that $$e^{\alpha_1 J + \alpha_2 \eta} = \alpha_2 \kappa +1,$$ for example $\alpha_2=1/\kappa$ and $\alpha_1= \frac{1}{J}\{ln(\alpha_2 \kappa +1)-\alpha_2 \eta\}$, and let $\alpha_3=\alpha_2$. Obviously, $\alpha_1 > 0$ when $\eta < (ln2)\kappa$.
Now, apply It\^o formula for the function $u(x,y,w)=e^{-\alpha_1(b-x)-\alpha_2(e-y)-\alpha_3(r-w)}$ to result 
\begin{align*}
\mathcal{L}_t^1u(Z_s, \lambda_s, D_{t,z}^{(2)}\lambda_s&)=-\alpha_2 \kappa \lambda_s u(Z_s,\lambda_s, D_{t,z}^{(2)}\lambda_s)-\alpha_3 \kappa D_{t,z}^{(2)}\lambda_s u(Z_s,\lambda_s, D_{t,z}^{(2)}\lambda_s)\\
&+u(Z_s,\lambda_s, D_{t,z}^{(2)}\lambda_s)\int_{\mathbb{R}_0} (e^{\alpha_1 \triangle Z_s+ \alpha_2 \triangle \lambda_s+ \alpha_3 \triangle D_{t,z}^{(2)}\lambda_s  }- 1)dz_0\\
&=-\alpha_2 \kappa \lambda_s u(Z_s,\lambda_s, D_{t,z}^{(2)}\lambda_s)-\alpha_3 \kappa D_{t,z}^{(2)}\lambda_s u(Z_s,\lambda_s, D_{t,z}^{(2)}\lambda_s)\\
&+ u(Z_s,\lambda_s, D_{t,z}^{(2)}\lambda_s)\Big[\int_{\mathbb{R}_0} \mathbb{I}_{(0,\lambda_s)}(z_0)(e^{\alpha_1 J+ \alpha_2 \eta }- 1)dz_0\\
& \        \     \qquad  +\int_{\mathbb{R}_0} \mathbb{I}_{(\lambda_s, \varphi_s(\bar{N}|_{(-\infty, s)}+\epsilon_{t,z}))}(z_0)(e^{\alpha_1 J+\alpha_3 \eta}- 1)dz_0\Big]\\
&= 0,
\end{align*}
and then, 
\begin{align}
&e^{-\alpha ({s \wedge\hat{\tau}_{t,b,e}})} u(Z_{s\wedge\hat{\tau}_{t,b,e}}, \lambda_{s \wedge\hat{\tau}_{t,b,e}}, D_{t,z}^{(2)}\lambda_{s \wedge\hat{\tau}_{t,b,e}}) -e^{-\alpha t} u(Z_{t}, \lambda_{t}, D_{t,z}^{(2)}\lambda_{t}) \nonumber\\
&  = \int_t^{s \wedge\hat{\tau}_{t,b,e}} e^{-\alpha v}\Big(-\alpha u(Z_v,\lambda_v, D_{t,z}^{(2)}\lambda_{v})+(\mu-\frac12 \sigma_1^2) \partial_1 u(Z_v,\lambda_v, D_{t,z}^{(2)}\lambda_{v})\Big) dv\nonumber\\
& + \int_t^{s \wedge\hat{\tau}_{t,b,e}} e^{-\alpha v} \frac12  \sigma_1^2 \frac{\partial^2}{\partial z^2}  u(Z_v,\lambda_v, D_{t,z}^{(2)}\lambda_{v})dv\nonumber\\
&+  \int_t^{s \wedge\hat{\tau}_{t,b,e}}  e^{-\alpha v}\Big(\kappa\Theta \partial_2 u(Z_v,\lambda_v, D_{t,z}^{(2)}\lambda_{v}) +  \mathcal{L}_t^1 u(Z_v,\lambda_v, D_{t,z}^{(2)}\lambda_{v})\Big) dv \nonumber\\
&+  \mathcal{M}^{(1)}_{s \wedge\hat{\tau}_{t,b,e}}+\mathcal{M}^{(2)}_{s \wedge\hat{\tau}_{t,b,e}},
\label{expectu}
\end{align}
where $\mathcal{M}^{(1)}$ nad $\mathcal{M}^{(2)}$ are continuous and non-continuous martingales. 
Take expectation on bothsides of above equality and let 
$$\alpha= (\mu-\frac12 \sigma_1^2) \alpha_1 +\frac12 \sigma_1^2\alpha_1^2 +\kappa \Theta \alpha_2,$$
 which will vanish the drift term of Equation (\ref{expectu}), and then tend $s$ to infinity to result
\begin{align*}
\mathbb{E}(e^{-\alpha \hat{\tau}_{t,b,e}}\vert \mathcal{F}_t)&= \mathbb{E}(e^{-\alpha {\hat{\tau}_{t,b,e}}} u(Z_{\hat{\tau}_{t,b,e}}, \lambda_{\hat{\tau}_{t,b,e}}, D_{t,z}^{(2)}\lambda_{\hat{\tau}_{t,b,e}}) \vert \mathcal{F}_t)\\
& =e^{-\alpha t} u(Z_{t}, \lambda_{t}, \eta \mathbb{I}_{(0, \lambda_t)}(z)).
\end{align*}
Now, the above expression and the following Laplace transform, which is convenient for numerical Laplace inversion, 
\begin{align*}
\frac{1}{\alpha} \mathbb{E} (e^{-\alpha \hat{\tau}_{t,b,e}}\vert \mathcal{F}_t) &= \int_t^\infty e^{-\alpha u} P(\hat{\tau}_{t,b,e} \leq u \vert \mathcal{F}_t) du \\
&= e^{-\alpha t} \int_0^\infty e^{-\alpha u} P(\hat{\tau}_{t,b,e}-t \leq u \vert \mathcal{F}_t) du,
\end{align*}   
will deduce the probability distribution of the processes $Z$ and $\lambda$ and $D_{t,z}\lambda$. These facts result
\begin{align*}
P(&\sup_{t \leq s \leq T}Z_s \geq y  \vert \mathcal{F}_t)\\
&= \int_0^\infty \int_0^\infty  P(\sup_{t \leq s \leq T}Z_s \geq y , \sup_{t \leq s \leq T} \lambda_s \geq e, \sup_{t \leq s \leq T} D_{t,z}\lambda_s \geq r\vert \mathcal{F}_t) dedr \\
&=\int_0^\infty  \int_0^\infty  P(\hat{\tau}_{t,y,e} \leq T \vert \mathcal{F}_t ) de dr\\
&=   L^{-1}(\frac{1}{u})(T-t)    \int_0^\infty \int_0^\infty e^{-\alpha_1 (y-Z_t) -\alpha_2(e-\lambda_t)-\alpha_2(r-D_{t,z}\lambda_t)} dedr\\
&=  \frac{1}{\alpha_2^2} L^{-1}(\frac{1}{u})(T-t)  e^{-\alpha_1 (y-Z_t) +\alpha_2\lambda_t+\alpha_2\eta\mathbb{1}_{(0,\lambda_t]}(z)},
\end{align*}
which $L^{-1}(f)(x)$ is the inverse Laplace operator of the function $f$ in the point $x$.
Take the expectation and use the fact that $\alpha_1, \alpha_2 >0$ to derive 
\begin{align}
\int_{M_t-K_{t,z}}^\infty P(&\sup_{t \leq s \leq T}Z_s \geq y \vert \mathcal{F}_t)dy\nonumber\\
& = \frac{1}{\alpha_2^2} L^{-1}(\frac{1}{u})(T-t)  e^{\alpha_2\lambda_t+\alpha_2\eta\mathbb{1}_{(0,\lambda_t]}(z)}   \int_{M_t-K_{t,z}}^\infty  e^{-\alpha_1 (y-X_t)} dy\nonumber \\
&=  \frac{1}{\alpha_2^2\alpha_1} L^{-1}(\frac{1}{u})(T-t)  e^{\alpha_2\lambda_t+\alpha_2\eta\mathbb{I}_{(0,\lambda_t]}(z)}    e^{-\alpha_1 (M_t-K_{t,z}-X_t)}.\label{distribute1}
\end{align}
In the rest of this subsection, we will present the distribution of the supremum of $X$ in the same manner. \\
With the same contribution, one can show that 
\begin{align*}
\mathcal{L}_t^1&u(X_s, \lambda_s) \\
&= -\alpha_2 \kappa \lambda_s u(X_s,\lambda_s)+ u(X_s,\lambda_s)\int_{\mathbb{R}_0} \mathbb{I}_{(0,\lambda_s)}(z_0)(e^{\alpha_1 J+ \alpha_2 \eta }- 1)dz_0=0,
\end{align*}
and then, for $\alpha= (\mu-\frac12 \sigma_1^2) \alpha_1 +\frac12 \sigma_1^2\alpha_1^2 +\kappa \Theta \alpha_2,$ using the It\^o formula for the function $u(x,y)=e^{-\alpha_1(b-x)-\alpha_2(e-y)}$ we drive
\begin{align}
e^{-\alpha ({s\wedge{\tau}_{t}(b,e)})} & u(X_{s\wedge{\tau}_{t}(b,e)}, \lambda_{s \wedge{\tau}_{t}(b,e)}) -e^{-\alpha t} u(X_{t}, \lambda_{t}) \nonumber\\
&  = \int_t^{s \wedge{\tau}_{t}(b,e)} e^{-\alpha v}\Big(\{-\alpha+(\mu-\frac12 \sigma_1^2)\alpha_1+ \frac12  \sigma_1^2\alpha_1^2\} u(X_v,\lambda_v)\Big) dv\nonumber\\
&+  \int_t^{s \wedge{\tau}_{t}(b,e)}  e^{-\alpha v}\Big(\kappa\Theta \alpha_2 u(X_v,\lambda_v) +  \mathcal{L}_t^1 u(X_v,\lambda_v)\Big) dv +\mathcal{M}_{t,s}\nonumber\\
&= \mathcal{M}_{t,s},\nonumber\\
\end{align}
which $\mathcal{M}_{t,s}$ is a martingale. 
Therefore, 
\begin{align*}
\mathbb{E}(e^{-\alpha {\tau}_{t}(b,e)}\vert \mathcal{F}_t)= \mathbb{E}(e^{-\alpha {{\tau}_{t}(b,e)}} u(X_{{\tau}_{t}(b,e)}, \lambda_{{\tau}_{t}(b,e)}) \vert \mathcal{F}_t) =e^{-\alpha t} u(X_{t}, \lambda_{t}),
\end{align*}
and using Laplace inversion we deduce
\begin{align*}
P(\sup_{t \leq s \leq T}X_s \geq b \vert \mathcal{F}_t)& = \int_0^\infty  P(\sup_{t \leq s \leq T}X_s \geq b , \sup_{t \leq s \leq T} \lambda_s \geq e\vert \mathcal{F}_t) de \\
&=\int_0^\infty  P({\tau}_{t}(b,e) \leq T \vert \mathcal{F}_t ) de \\
&=   L^{-1}(\frac{1}{u})(T-t)   \int_0^\infty e^{-\alpha_1 (b-X_t) -\alpha_2(e-\lambda_t)} de\\
&=  \frac{1}{\alpha_2} L^{-1}(\frac{1}{u})(T-t)  e^{-\alpha_1 (b-X_t) +\alpha_2\lambda_t},
\end{align*}
which results
\begin{align}
\int_{M_t}^\infty P(\sup_{t \leq s \leq T}X_s \geq x \vert \mathcal{F}_t)dx& = \frac{1}{\alpha_2} L^{-1}(\frac{1}{u})(T-t)  e^{\alpha_2\lambda_t}   \int_{M_t}^\infty  e^{-\alpha_1 (x-X_t)} dx \nonumber\\
&=  \frac{1}{\alpha_2\alpha_1} L^{-1}(\frac{1}{u})(T-t)  e^{\alpha_2\lambda_t}    e^{-\alpha_1 (M_t-X_t)}. \label{distribute2}
\end{align}
Substitute \eqref{distribute1} and \eqref{distribute2} in \eqref{distribute} to have
\begin{align*}
\mathbb{E}\left[D_{t, z}^{(2)} F \mid \mathcal{F}_{t}\right]= \frac{1}{\alpha_2\alpha_1} L^{-1}(\frac{1}{u})(T-t)  e^{\alpha_2\lambda_t}  e^{-\alpha_1 (M_t-X_t)} [e^{\alpha_1 K_{t,z}}-1].
\end{align*}




\section{Appendix}
We will prove the Lemma \ref{distcox1}. With the same contribution of the previous section, in the Cox process,  one can show that for every $t \leq s$
\begin{align*}
\mathcal{L}_t^0 u(X_s, \lambda_s)
&:= \Big[ -\alpha_2 \kappa +\frac12 \sigma_2^2 \alpha_2^2-\alpha_1 (\bar{\mu}-J)\Big]\lambda_s u(X_s,\lambda_s)\\
&+ u(X_s,\lambda_s)\int_{\mathbb{R}_0} \Big[ u(X_s+\triangle X_s, \lambda_s)-u(X_s,\lambda_s)\Big]\nu(dz).
\end{align*}
Set $\alpha_2 = \frac{2\kappa (\Theta +1)}{\sigma_2^2}$, that it results $ -\alpha_2 \kappa +\frac12 \sigma_2^2 \alpha_2^2 =\alpha_2 \kappa \Theta$, and also choose some positive constant more than one for $\alpha_1$ and then let  $$\alpha= (\mu-\frac12 \sigma_1^2) \alpha_1 +\frac12 \sigma_1^2\alpha_1^2 +\alpha_1 (\bar{\mu}-J)+e^{\alpha_1 J}-1.$$ 
Applying the It\^o formula for the function $u(x,y)=e^{-\alpha_1(b-x)-\alpha_2(e-y)}$ we drive
\begin{align}
e^{-\alpha ({s\wedge{\tau}_{t}^0 (b,e)})} & u(X_{s\wedge{\tau}_{t}^0(b,e)}, \lambda_{s \wedge{\tau}_{t}^0(b,e)}) -e^{-\alpha t} u(X_{t}, \lambda_{t}) \nonumber\\
&  = \int_t^{s \wedge{\tau}_{t}^0(b,e)} e^{-\alpha v}\Big(\{-\alpha+(\mu-\frac12 \sigma_1^2)\alpha_1+ \frac12  \sigma_1^2\alpha_1^2\} u(X_v,\lambda_v)\Big) dv\nonumber\\
&+  \int_t^{s \wedge{\tau}_{t}^0(b,e)}  e^{-\alpha v}\Big(\alpha_2 \kappa \Theta+ \mathcal{L}_t^0 u(X_v,\lambda_v)\Big) dv +\mathcal{M}_{t,s}\nonumber\\
&= \mathcal{M}_{t,s},\nonumber\\
\end{align}
which $\mathcal{M}_{t,s}$ is a martingale and ${\tau}_{t}^0 (b,e)= \inf\{t \leq s \leq T;  X_s \leq b~ or~ \lambda_s \leq e\}$. 
Therefore, 
\begin{align*}
\mathbb{E}(e^{-\alpha {\tau}_{t}^0(b,e)}\vert \mathcal{F}_t)= \mathbb{E}(e^{-\alpha {{\tau}_{t}^0(b,e)}} u(X_{{\tau}_{t}^0(b,e)}, \lambda_{{\tau}_{t}(b,e)}) \vert \mathcal{F}_t) =e^{-\alpha t} u(X_{t}, \lambda_{t}),
\end{align*}
and using Laplace inversion we deduce
$$\bar{F}_{t,x,y} =   P({\tau}_{t}^0(x,y) \leq T \vert \mathcal{F}_t ) =  \frac{1}{\alpha_2} L^{-1}(\frac{1}{u})(T-t)  e^{-\alpha_1 (x-X_t)  -\alpha_2(y-\lambda_t)}.$$
%
%
\\

\end{document}